\documentclass[a4paper,9pt]{article}

\usepackage{amsmath,amscd,amssymb,epsfig,xypic}

\def\g{\mathfrak{g}}
\def\p{\mathfrak{p}}
\def\q{\mathfrak{q}}
\def\a{\mathfrak{a}}
\def\b{\mathfrak{b}}
\def\ph{\varphi}
\def\oa{\overline{\alpha}}
\def\ob{\overline{\beta}}
\def\oj{\overline{\gamma}}
\def\e{\frac{1}{2}}
\def\h{\mathfrak{h}}
\def\n{\mathfrak{n}}

\def\l{\mathfrak{l}}

\def\C{\mathbb{C}}
\def\Z{\mathbb{Z}}

\newtheorem{theo}{\bf{Theorem}}
\newtheorem{prop}{Proposition}
\newtheorem{lem}{Lemma}
\newtheorem{cor}{Corollary}

\begin{document}
\title{On classification of dynamical r-matrices}
\author{Olivier Schiffmann\\
Harvard University and ENS Paris\\
45 rue d'Ulm,\\
75 005 Paris\\
\texttt{schiffma@clipper.ens.fr}}
\maketitle
\begin{abstract}
Using the  gauge transformations of the Classical Dynamical Yang-Baxter
Equation
introduced by P. Etingof and A. Varchenko in \cite{EV}, we reduce the
classification of dynamical
r-matrices $r$ on a commutative subalgebra $\l$ of a Lie algebra $\g$ to a
purely algebraic
problem, under some assumption on the symmetric part of $r$. We then describe, for a simple complex Lie algebra $\g$, all non
skew-symmetric
dynamical r-matrices on a commutative subalgebra  $\l \subset\g$ which
contains a regular
semisimple element. This interpolates results of P. Etingof and A.
Varchenko (\cite{EV}, when $\l$ is a
Cartan subalgebra) and results of A. Belavin and V. Drinfeld for constant
r-matrices (\cite{BD}).
This classification is similar, and in some sense simpler than the
Belavin-Drinfeld
classification.
\end{abstract}
\section{The Classical Yang-Baxter Equation}
\paragraph{}Let $\g$ be a Lie algebra. The CYBE is the following algebraic
equation for an element $r \in \g \otimes \g$:
\begin{equation}\label{E:01}
[r^{12},r^{13}]+[r^{12},r^{23}]+[r^{13},r^{23}]=0.
\end{equation}
Solutions of this   equation are called r-matrices. In the theory of
quantum groups,
one is mainly interested in r-matrices satisfying
\begin{equation}\label{E:02}
r+r^{21} \in (S^2 \g)^{\g}.
\end{equation}
  See \cite{CP} for the links with the theory of quantum groups, and
\cite{C} for
links with Conformal Field Theory and the Wess-Zumino-Witten model on
$\mathbb{P}^1$. The
geometric interpretation of the CYBE was given by Drinfeld in terms of
Poisson-Lie groups
(\cite{Dr1}).
\section{The Belavin-Drinfeld Classification}
\paragraph{Notations:}Let $\g$ be a simple complex Lie algebra with a
nondegenerate invariant form $(\,,\,)$, $\h \subset \g$ a Cartan subalgebra
and $\Delta$ the root system. For $\alpha \in \Delta$, let $\g_{\alpha}$
denote the root subspace associated to $\alpha$. Let $W$ be the Weyl group
and $s_{\alpha}$, $\alpha \in \Delta$ the reflection with respect to
$\alpha^{\perp}$. Finally, let $\Omega\in S^{2}\g$ and $\Omega_{\h}\in
S^2\h$ be the inverse elements to the form $(\,,\,)$. Notice that
$(S^2\g)^{\g} = \C \Omega$.\\
\hbox to1em{\hfill}For any polarization $\g=\n_{-} \oplus \h \oplus
\n_{+}$, we denote by $\Pi$ or $\Pi(\n_{+})$ the corresponding set of
simple positive roots, by $\Delta_+$ the set of positive roots and by
$\b_{\pm}=\n_{\pm} \oplus \h$ the Borel subalgebras. For $\Gamma \subset
\Pi$, set $\left<\Gamma\right>=\Z \Gamma \cap \Delta$, and let
$\g_{\Gamma}$ be the subalgebra generated by $\g_{\alpha}$, $\alpha \in
\left<\Gamma\right>$. We will write $\g_{\Gamma}=\n_{+}(\Gamma) \oplus
\h(\Gamma) \oplus \n_{-}(\Gamma)$ for the induced polarization and
$W(\Gamma)$ for the subgroup of $W$ generated by $s_{\alpha}$, $\alpha \in
\Gamma$.
\paragraph{}Let us fix a polarization of $\g$.\\
{\bf{Definition:}} A {\it{Belavin-Drinfeld triple}} is a triple
$(\Gamma_{1},\Gamma_{2},\tau)$ where $\Gamma_{1},\Gamma_{2} \subset \Pi$
and $\tau: \Gamma_{1} \stackrel{\sim}{\to} \Gamma_{2}$ is a norm-preserving
bijection satisfying the following "nilpotency" condition:\\
\hbox to1em{\hfill}"For any $\gamma_{1} \in \Gamma_{1}$, there exists $n
>0$ such that $\tau^{n}(\gamma_{1}) \in \Gamma_{2} \backslash \Gamma_{1}$".
\paragraph{}Let $(\Gamma_{1},\Gamma_{2},\tau)$ be a Belavin-Drinfeld
triple. For each choice of Chevalley generators
$(e_{\alpha},f_{\alpha},h_{\alpha})_{\alpha \in \Gamma_{i}}$, $i=1,2$, the
isomorphism $\tau$ induces a Lie algebra isomorphism $\g_{\Gamma_{1}}
\stackrel{\sim}{\to} \g_{\Gamma_{2}}$ (by $e_{\alpha} \mapsto
e_{\tau(\alpha)},\,f_{\alpha} \mapsto f_{\tau(\alpha)},\,h_{\alpha} \mapsto
h_{\tau(\alpha)}$).\\
 Define a partial order on $\Delta_{+}$ by setting $\alpha < \beta$ if
there exists $n >0$ such that $\tau^n(\alpha)=\beta$ (in particular,
$\alpha \in \Gamma_{1}$ and $\beta \in \Gamma_{2}$).
\paragraph{Definition:} A basis $(x_{\alpha})_{\alpha \in \Delta}$ of
$\n_{+} \oplus \n_{-}$ is called {\it{admissible}} if
$(x_{\alpha},x_{-\alpha})=1$ and $\tau(x_{\alpha})=x_{\tau(\alpha)}$ for
$\alpha \in \left<\Gamma_{1}\right>$.

\begin{theo}[Belavin-Drinfeld] Let $\g$ be a simple complex Lie algebra.
\begin{enumerate}
\item Let $(\Gamma_{1},\Gamma_{2},\tau)$ be a Belavin-Drinfeld triple,
$(x_{\alpha})$ an admissible basis, and let $r_{0} \in \h \otimes \h$ be
such that
\begin{equation}\label{E:03}
r_{0} + r_{0}^{21} = \Omega_{\h},
\end{equation}
\begin{equation}\label{E:04}
(\tau(\alpha) \otimes 1)r + (1 \otimes \alpha)r =0\qquad
\mathrm{for}\;\alpha \in \Gamma_{1}.
\end{equation}
Then
\begin{equation}\label{E:05}
r=r_{0} + \sum_{\alpha \in \Delta_{+}} x_{-\alpha} \otimes x_{\alpha} +
\sum_{\alpha,\beta \in \Delta_{+}, \alpha < \beta} x_{-\alpha} \wedge
x_{\beta}
\end{equation}
is an r-matrix satisfying $r + r^{21} = \Omega$.
\item Any r-matrix satisfying $r + r^{21}=\Omega$ is of the above type for
a suitable polarization of $\g$.
\end{enumerate}
\end{theo}
This theorem is proved in \cite{BD}. For instance, the standard r-matrix
for a fixed polarization $r=\frac{\Omega_{\h}}{2}+ \sum_{\alpha \in
\Delta_+}
x_{-\alpha} \otimes x_{\alpha}$ corresponds to
$\Gamma_{1}=\Gamma_{2}=\emptyset$.\\
{\bf{Remark:}} Skew-symmetric r-matrices admit a well known interpretation
in terms of
nondegenerate 2-cocycles on Lie subalgebras of $\g$ (\cite{Dr1}), but their
classification is
unavailable since it requires a classification of Lie subalgebras in $\g$.
\section{The Dynamical Yang-Baxter Equation}
\paragraph{}Let $\g$ be a Lie algebra over $\C$ and $\l \subset \g$ a
subalgebra. An
element $x \in \g \otimes \g$ will be called $\l$-invariant if
\begin{equation}\label{E:07}
[k\otimes 1 + 1 \otimes k,x]=0 \qquad \;\;\;\;(\forall\; k \in \l).
\end{equation}
For $x \in \g^{\otimes 3}$, we let $\mathrm{Alt}(x)=x^{123} + x^{231}
+x^{312}$. Let $D
\subset
\l^*$ be any open region.\\
\hbox to1em{\hfill}The CDYBE is the following differential equation for a
holomorphic
$\l$-invariant function
$r:\;\ D \to \g \otimes
\g$:
\begin{equation}\label{E:06}
\mathrm{Alt}(dr) + [r^{12},r^{13}]+[r^{12},r^{23}]+[r^{13},r^{23}]=0,
\end{equation}
where the differential of $r$ is considered as a holomorphic function
$$dr: D \to \g \otimes \g \otimes \g,\qquad \lambda \mapsto \sum_{i} x_{i}
\otimes
\frac{\partial r^{23}}{\partial x_{i}}(\lambda),\qquad (\lambda \in
\l^*),$$ for any basis
$(x_{i})$ of $\l$. In this case,
$$\mathrm{Alt}(dr)=\sum_{i}x_{i}^{(1)} \frac{\partial r^{23}}{\partial x_{i}} +
\sum_{i}x_{i}^{(2)} \frac{\partial r^{31}}{\partial x_{i}} +
\sum_{i}x_{i}^{(3)} \frac{\partial r^{12}}{\partial x_{i}}.$$
The solutions to this equation are called {\it{dynamical}} r-matrices.
Dynamical r-matrices which are relevant to the theory of quantum groups are
those satisfying the following condition, analogous to~(\ref{E:02}):
\begin{equation}
\mathrm{Generalized\; unitarity:}\; r(\lambda) + r^{21}(\lambda) \in
(S^2\g)^{\g}. \label{E:08}
\end{equation}
{\bf{Remark:}} the CDYBE was first written down by G. Felder and C.
Wiezcerkowski in connection with the Wess-Zumino-Witten model on elliptic
curves (\cite{FW}). The relation with elliptic quantum groups is explained
in \cite{Fe}. A geometric interpretation of the CDYBE analogous to the
theory of Poisson-Lie groups for the CYBE is given in \cite{EV}.
\section{Gauge transformations:}
\paragraph{}We recall some results from \cite{EV}. We suppose here that
$\l$ is commutative
and we let $D$ be the formal polydisc centered at the origin. Let
$G$ be a complex Lie group such that
$\mathrm{Lie}(G)=\g$, and let $L$ be the connected subgroup of $G$ such that
$\mathrm{Lie}(L)=\l$. Let $G^L$ be the centralizer of $L$ in $G$ and
$\g^{\l}$ its Lie algebra.We will denote by $(\g \otimes \g)^{\l}$ the space of all $\l$-invariant elements in $\g \otimes \g$.\\
\hbox to1em{\hfill}Let $g: D \to G^L$ be any holomorphic function; the $1$-form
$\eta=g^{-1}dg$ gives rise to a function $ \overline{\eta}: D \to \l
\otimes \g^{\l}$. If $r: D \to
(\g\otimes
\g)^{\l}$ is an $\l$-invariant function satisfying~(\ref{E:08}), we set
$$r^g=(g \otimes g)(r-\overline{\eta}+\overline{\eta}^{21})(g^{-1} \otimes
g^{-1}).$$
\begin{prop} The function $r$ is a dynamical r-matrix if and only if the
function $r^g$ is.
\end{prop}
\hbox to1em{\hfill}Thus the group $\mathrm{Map}(D,G^L)$ is a gauge
transformation group for the CDYBE. Notice that this group is not
commutative if $G^L$ isn't.
\begin{theo}\label{T:EV} Let $\rho,r:D \to \g^{\otimes 2}$ be two dynamical
r-matrices satisfying~(\ref{E:08}) such that $r(0)=\rho(0)$. Then there
exists $g \in \mathrm{Map}(D,G^L)$ such that $\rho=r^g$.
\end{theo}
\hbox to1em{\hfill}This shows that the space of dynamical r-matrices is, up
to gauge equivalence, finite dimensional. Proofs of the above results can
be found in \cite{EV}.
\paragraph{}We will now prove a converse of Theorem~\ref{T:EV} which
reduces the CDYBE to a purely algebraic
equation under some assumption on the symmetric part $\frac{\Omega}{2}$ of $r$: let $\Omega \in (S^2\g)^{\g}$, let $\g_{\Omega}$ be the ideal in $\g$ generated by the components of $\Omega$ and denote by $\g_{\Omega}=\bigoplus_{\lambda} \g_{\Omega}(\lambda)$ the generalized weight space decomposition of $\g_{\Omega}$ with respect to the adjoint action of $\l$. The condition we will need is the following:
\begin{equation}
\g^\l\;\mathrm{acts\;semisimply\;on\;} \g_{\Omega}(0) \tag{*}
\end{equation}
\hbox to1em{\hfill}Suppose that (*) is fulfilled and let $z(\g^\l)$ denote the center of $\g^\l$. Then we have a decomposition $\g_{\Omega}(0)=z_0(\g^\l) \oplus V$ where $z_0(\g^\l)=z(\g^\l) \cap \g_\Omega (0)$ and $V$ is the sum of all non-trivial irreducible $\g^\l$-modules in $\g_\Omega(0)$. It is clear that $\l \cap V = \{0\}$. We will say that a complement $\l'$ of $\l$ in $\g$ is admissible if $V \subset \l'$, and write $\pi: \g \to \l$ for the projection along $\l'$. Notice that by $\g^\l$-invariance of $\Omega$, 
\begin{equation}\label{E:decompomega}
\Omega \in S^2 z_0(\g^\l) \oplus S^2 V \oplus \bigoplus_{\lambda \neq 0} \g_{\Omega}(\lambda) \otimes \g_{\Omega}(-\lambda).
\end{equation}
We will denote by
$CYB: \g^{\otimes 2} \to \g^{\otimes 3}$ the map:
$$r \mapsto  [r^{12},r^{13}]+[r^{12},r^{23}]+[r^{13},r^{23}].$$
\hbox to1em{\hfill}It is more convenient to work with the skew-symmetric
part of $r$. If $r(\lambda) + r^{21}(\lambda)=\Omega \in (S^2(\g))^{\g}$,
we set $s(\lambda)=r(\lambda)-\frac{\Omega}{2}$. It is easy to see that the
CDYBE for $r$ is equivalent to the following equation for $s$:
\begin{equation}\label{E:g02}
\mathrm{Alt}(ds) + CYB(s) + \frac{1}{4}CYB(\Omega)=0.
\end{equation}
Recall that as $\Omega$ is symmetric and invariant,
$CYB(\Omega)=[\Omega_{13},\Omega_{23}]$.\\
\begin{theo}\label{T:g01}
Let $G$ be a complex Lie group and $L \subset G$ a connected commutative
subgroup. Let
$\g,\l,\g^{\l}$ denote the Lie algebras of $G,L$ and $G^L$. Let
$\Omega \in (S^2 \g)^{\g}$.
Then
\begin{enumerate}
\item Let $\l'$ be any complement of $\l$ in $\g$. Any dynamical r-matrix $r(\lambda)$ on $\l$ such that $r(\lambda)
+r^{21}(\lambda)=\Omega$ is gauge equivalent to a dynamical r-matrix
$\tilde{r}(\lambda): D \to \frac{\Omega}{2} +
(\Lambda^{2}(\l'))^{\l}$.
\item Suppose that condition (*) is true and let $\l'$ be any admissible complement of $\l$ in $\g$. Let $r_0 \in \frac{\Omega}{2} + (\Lambda^{2}(\l'))^{\l}$ satisfy
 \begin{equation}\label{E:g01}
CYB(r_0) \in \mathrm{Alt}(\l \otimes \g \otimes \g)
\end{equation}
such that $s_0=r_0-\frac{\Omega}{2}$ is a regular point of the algebraic manifold
$$M_{\Omega}=\{ s \in (\Lambda^2(\l'))^{\l}\,|\,CYB(s+ \frac{\Omega}{2})
\in \mathrm{Alt}(\l \otimes \g \otimes \g)\}.$$
Then there exists a dynamical r-matrix $r(\lambda): D \to  \frac{\Omega}{2} + (\Lambda^2(\l'))^{\l}$ such that $r(0)=r_0$.
\end{enumerate}
\end{theo}
\hbox to1em{\hfill} The condition (*) is satisfied in the following two interesting special cases: when $\Omega=0$ (triangular case) or when $\g^\l$ acts semisimply on $\g$ (for instance, $G$ is reductive and $L$ is contained in a maximal torus of $G$ or more generally, if $G^L$ is reductive).\\ 
\hbox to1em{\hfill}The proof of this theorem will occupy the rest of this
section.\\
Let us first prove part 1:
\begin{lem} Any dynamical r-matrix such that $r(\lambda) +
r^{21}(\lambda)=\Omega$ is gauge-equivalent to a dynamical r-matrix
$\tilde{r}(\lambda)$ such that $\tilde{r}(0) \in \frac{\Omega}{2} +
(\Lambda^2(\l')^{\l}$.
\end{lem}
{\bf{Proof:}} Let $\overline{\eta} \in \l \otimes \g^{\l}$ be such that
$r(0)-\overline{\eta} +
\overline{\eta}^{21}
\in
\frac{\Omega}{2} + \Lambda^2(\l')$. There exists a function $g: D \to G^L$
such that
$g^{-1}dg (0) = \eta$ (see \cite{EV}, Lemma 1.3). It is easy to see that
$\tilde{r}=r^g$ satisfies
the desired conditions.
\begin{flushright}
$\square\;\;\;\;\;\;\;\;$
\end{flushright}
\paragraph{}By Theorem~\ref{T:EV} , part 1. is proved. Let us now prove part 2. We will interpret the CDYBE (\ref{E:g02}) as a consistent system of differential equations defined on $M_{\Omega}$.\\
\hbox to1em{\hfill}For $s \in M_{\Omega}$, (\ref{E:g02}) is equivalent to
$$(\pi \otimes 1 \otimes 1) Alt(ds)=-(\pi \otimes 1 \otimes 1)(CYB(s) + \frac{1}{4}CYB (\Omega)).$$
This reduces to
\begin{equation}\label{E:g03}
ds=-(\pi \otimes 1 \otimes 1)([s^{12},s^{13}] + \frac{1}{4}CYB(\Omega)),
\end{equation}
or, in coordinates $(x_i)$, where $(x_i)$ is a basis of $\l$,
$$ \frac{\partial s}{\partial x_i}=-(x_i \otimes 1 \otimes 1)([s^{12},s^{13}] + \frac{1}{4}CYB(\Omega)).$$
\begin{lem}The system~(\ref{E:g03}) is consistent.
\end{lem}
{\bf{Proof:}} Set $X: M_{\Omega} \to \l \otimes \g \otimes \g,\;s \mapsto (\pi \otimes 1 \otimes 1)([s^{12},s^{13}]+ \frac{1}{4}CYB(\Omega))$. By definition, the curvature of (\ref{E:g03}) is given by
\begin{align*}
\sum_{i,j} &x_i \otimes x_j \otimes \big( \frac{\partial^2 s}{\partial x_i \partial x_j}-\frac{\partial^2 s}{\partial x_j \partial x_i} \big)\\
& = (\pi \otimes \pi \otimes 1 \otimes 1)\big(\big\{ [s^{23},[s^{12},s^{14}]] + [s^{23},\frac{1}{4}CYB(\Omega)^{124}]\\
&+ [[s^{12},s^{13}],s^{24}] + [\frac{1}{4}CYB(\Omega)^{123},s^{24}] \big\}\\
& - \big\{[s^{13},[s^{21},s^{24}]] + [s^{13},\frac{1}{4}CYB(\Omega)^{214}]\\
&+ [[s^{21},s^{23}],s^{14}] + [\frac{1}{4}CYB(\Omega)^{213},s^{14}]\big\} \big)\\
&= (\pi \otimes \pi \otimes 1 \otimes 1) \big( \big\{ [s^{23},[s^{12},s^{14}]] + [[s^{12},s^{13}],s^{24}] -[s^{13},[s^{21},s^{24}]] - [[s^{21},s^{23}],s^{14}]  \big\}\\
& +  \frac{1}{4} \big\{[s^{13} + s^{23}, CYB(\Omega)^{124}]-[s^{14} + s^{24}, CYB(\Omega)^{123}] \big\} \big).
\end{align*}
By the Jacobi identity, 
$$  [s^{23},[s^{12},s^{14}]] =[[s^{21},s^{23}],s^{14}], \qquad  [[s^{12},s^{13}],s^{24}] =[s^{13},[s^{21},s^{24}]]. $$
By $\g$-invariance of $CYB(\Omega)$, we have
\begin{align*}
[s^{13}+s^{23}, CYB(\Omega)^{124}]&=[s^{34}, CYB(\Omega)^{124}],\\
[s^{14}+s^{24}, CYB(\Omega)^{123}]&=-[s^{34}, CYB(\Omega)^{123}].
\end{align*}
Overall, we have the following expression for the curvature of (\ref{E:g03}):
$$\frac{1}{4}(\pi \otimes \pi \otimes 1 \otimes 1) ([CYB(\Omega)^{123}+CYB(\Omega)^{124}, s^{34}]=\frac{1}{4}[(\pi \otimes \pi \otimes 1)CYB(\Omega),s]$$
But (\ref{E:decompomega}) and the fact that $\l'$ is admissible imply that $(\pi \otimes \pi \otimes 1)CYB(\Omega)=0$. Thus, (\ref{E:g03}) is consistent.
\begin{flushright}
$\square\;\;\;\;\;\;\;\;$
\end{flushright}
\begin{lem}\label{L:g01} The system~(\ref{E:g03}) is defined on $M_{\Omega}$, i.e
the vector fields defined by~(\ref{E:g03}) are tangent to $M_{\Omega}$.
\end{lem}
{\bf{Proof:}} Let $x^* \in \l^* \stackrel{\pi^*}{\hookrightarrow} \g^*$, and set $h=(x^* \otimes 1 \otimes 1)([s^{12},s^{13}]+ \frac{1}{4}CYB(\Omega))$. Since $s \in
\Lambda^2(\l')$ we have $(x^* \otimes 1 \otimes 1)[s^{12},s^{13}] \in
\Lambda^2(\l')$.
Moreover, the admissibility of $\l'$ and (\ref{E:decompomega}) together imply that
$(x^* \otimes 1 \otimes 1)(CYB(\Omega))\in (\Lambda^2\l')^{\l}$
since \\
$[\l \otimes 1, S^2 z_0(\g^\l)]=0$. Thus $h \in \Lambda^2 \l'$.\\
\hbox to1em{\hfill}To conclude the proof of Lemma~\ref{L:g01} and
Theorem~\ref{T:g01}, we now show that
\begin{equation}\label{E:g04}
\begin{split}
[&s^{12},h^{13}] + [s^{12},h^{23}] + [s^{13},h^{23}]\\
 &+ [h^{12},s^{13}] + [h^{12},s^{23}] + [h^{13},s^{23}] \in
\mathrm{Alt}(\l \otimes \g \otimes \g).
\end{split}
\end{equation}
\hbox to1em{\hfill}To make the presentation more clear, we will use the
pictorial technique to represent expressions and make computations: we
associate to each morphism from a $n$-tensor to a $m$-tensor a diagram in
the following way: the operation of taking the commutator is represented by
$$
\centerline{\epsfbox{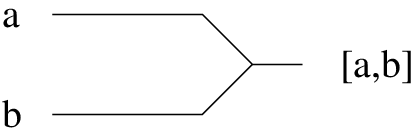}}
$$
Applying a linear form $x^*$ will be denoted by
$$
\centerline{\epsfbox{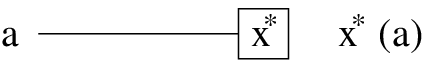}}
$$
Finally, we will represent $s$ and $\frac{\Omega}{2}$, which can be thought
of as maps from a $0$-tensor to a $2$-tensor, by
$$
\centerline{\epsfbox{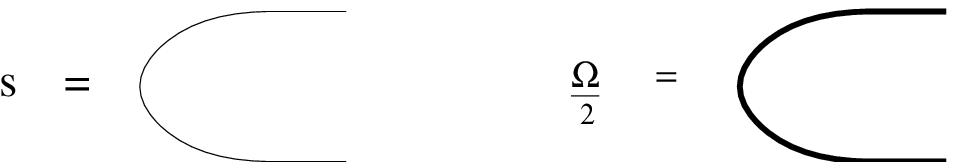}}
$$
For instance,
$$
\centerline{\epsfbox{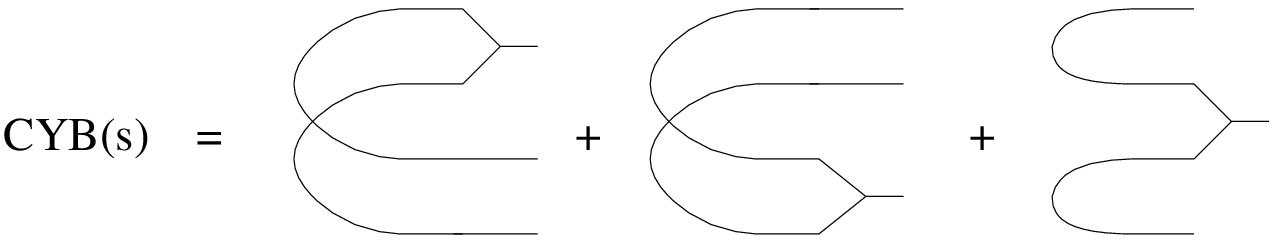}}
$$
\begin{lem}\label{L:g02} We have $x^{*(3)}
[CYB(s+\frac{\Omega}{2})^{123},s^{34}] \in \mathrm{Alt}(\l \otimes \g
\otimes \g)$ or, in pictures (modulo $\mathrm{Alt}(\l \otimes \g \otimes
\g)$)
$$
\centerline{\epsfbox{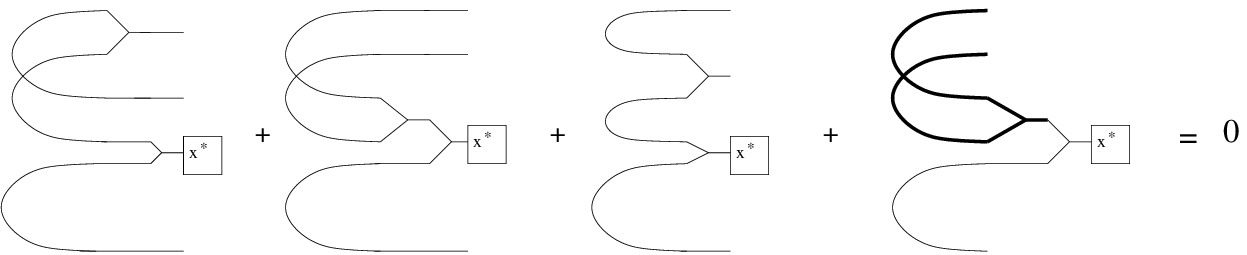}}
$$
\end{lem}
{\bf{Proof:}} Recall that $CYB(s+\frac{\Omega}{2}) \in \mathrm{Alt}(\l
\otimes \g \otimes \g)$. Thus the only part of the above expression which
can lie outside of $\mathrm{Alt}(\l \otimes \g \otimes \g)$ is obtained
from the $\g \otimes \g \otimes \l$-part of $CYB(s)$. But if $y \in \l$,
$$(x^* \otimes 1)[y \otimes 1, s]=-(x^* \otimes 1)[1 \otimes y,s]$$
by $\l$-invariance of $s$. This last expression is zero since $s \in
(\Lambda^2(\l'))^{\l}$. Lemma~\ref{L:g02} is proved.
\begin{flushright}
$\square \;\;\;\;\;\;\;\;$
\end{flushright}
It is clear how to generalize Lemma~\ref{L:g02} to other expressions of the
form
$$x^{*(k)}[CYB(s+\frac{\Omega}{2})^{123},s^{k4}].$$
\paragraph{}Now,~(\ref{E:g04}) can be drawn as
$$
\centerline{\epsfbox{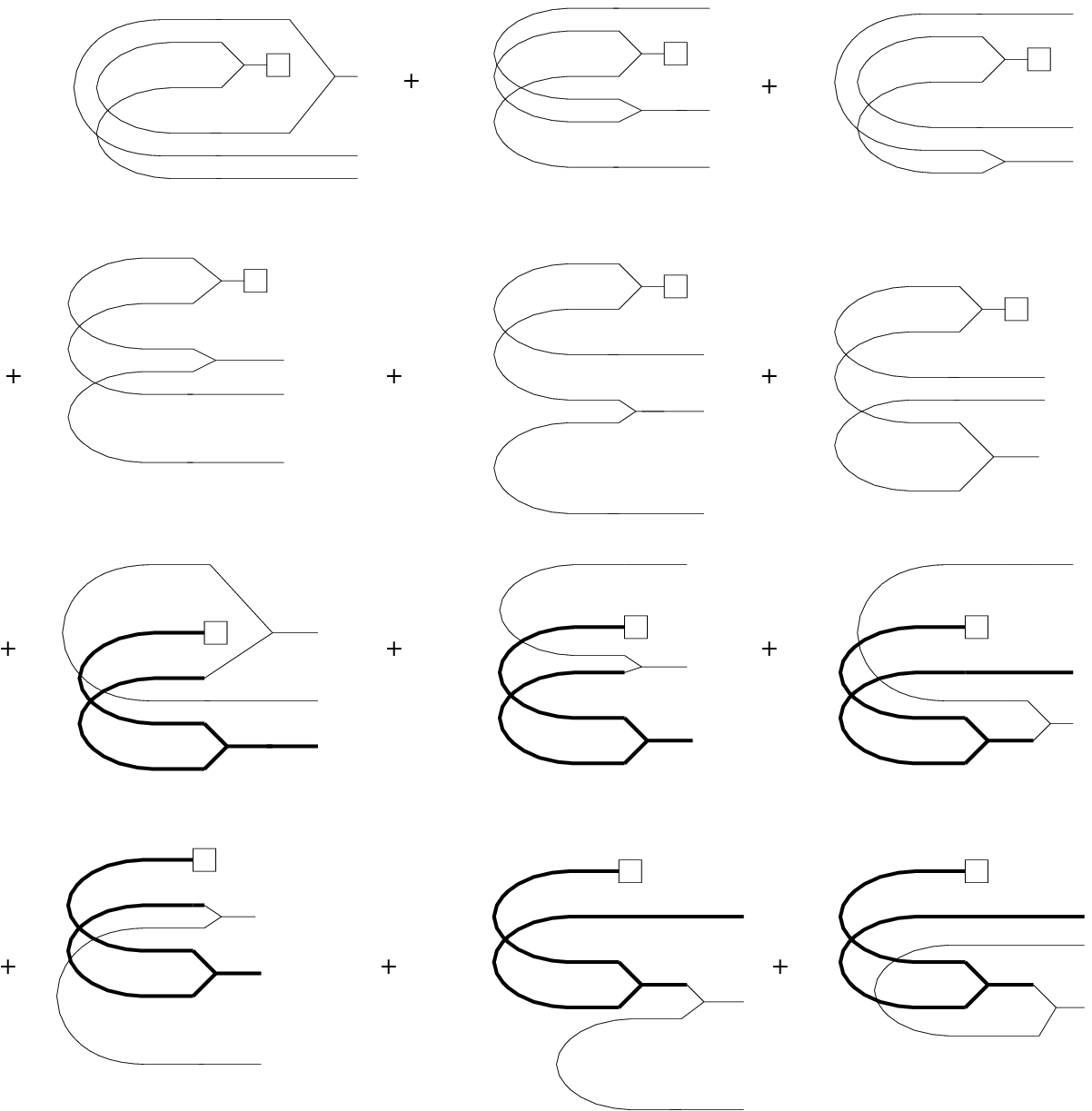}}
$$

but by Lemma~(\ref{L:g02}) we have, modulo $\mathrm{Alt}(\l \otimes \g
\otimes \g)$,
$$
\centerline{\epsfbox{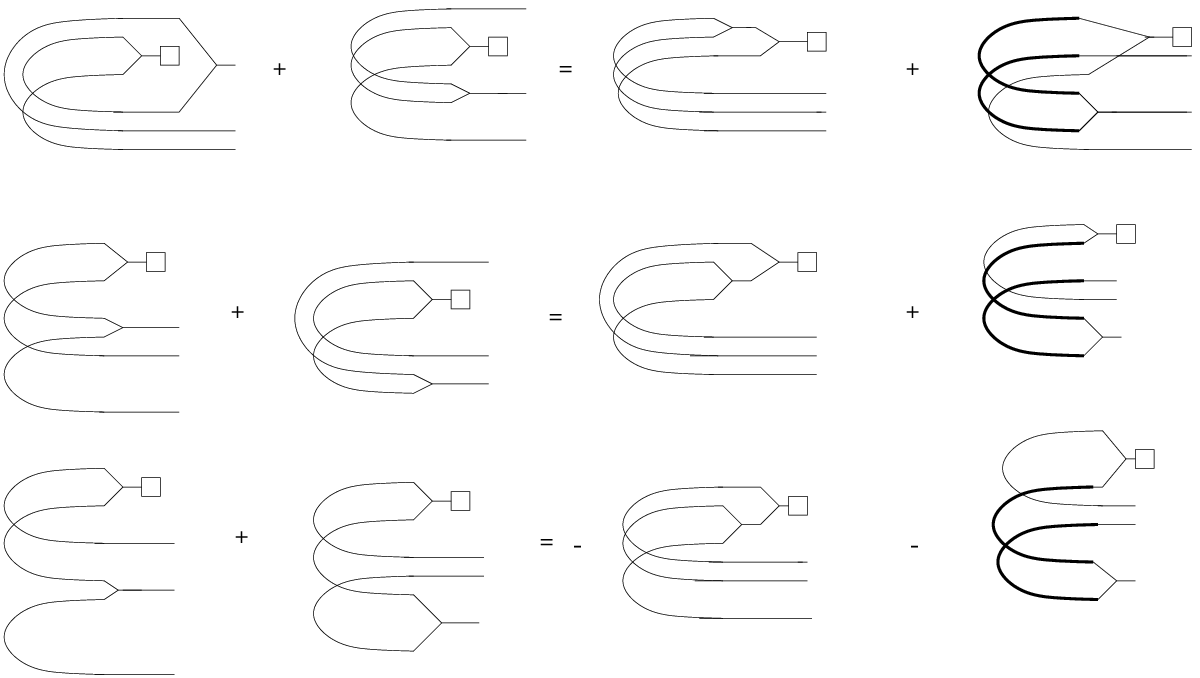}}
$$
It is easy to check that the sum of the terms of type $[CYB(s),s]$ in this
last expression is zero by the Jacobi identity. Moreover, by
$\g$-invariance of $\Omega$, we have
$$
\centerline{\epsfbox{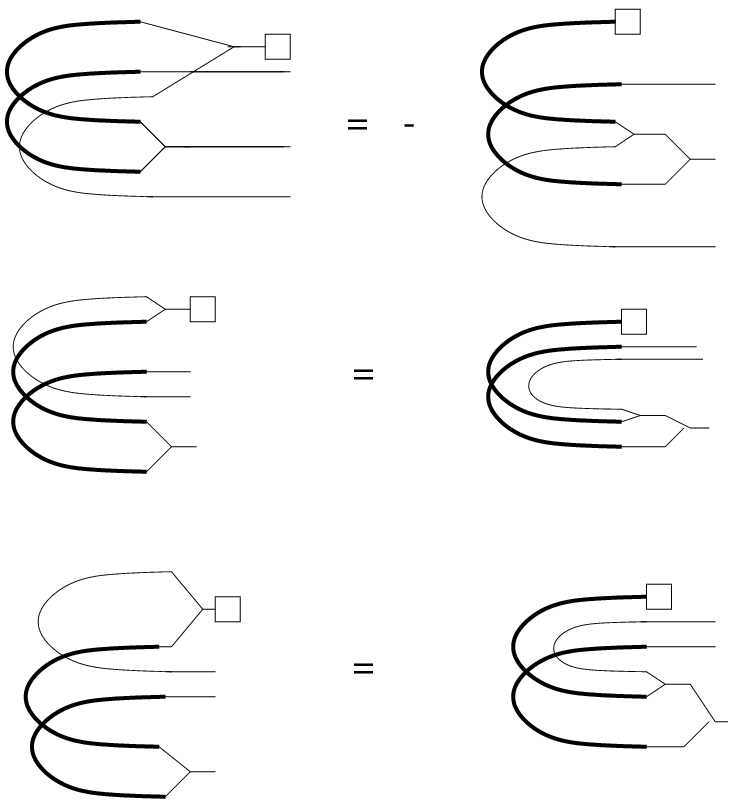}}
$$
Thus, modulo $\mathrm{Alt}(\l \otimes \g \otimes \g)$,~(\ref{E:g04}) reduces to
$$
\centerline{\epsfbox{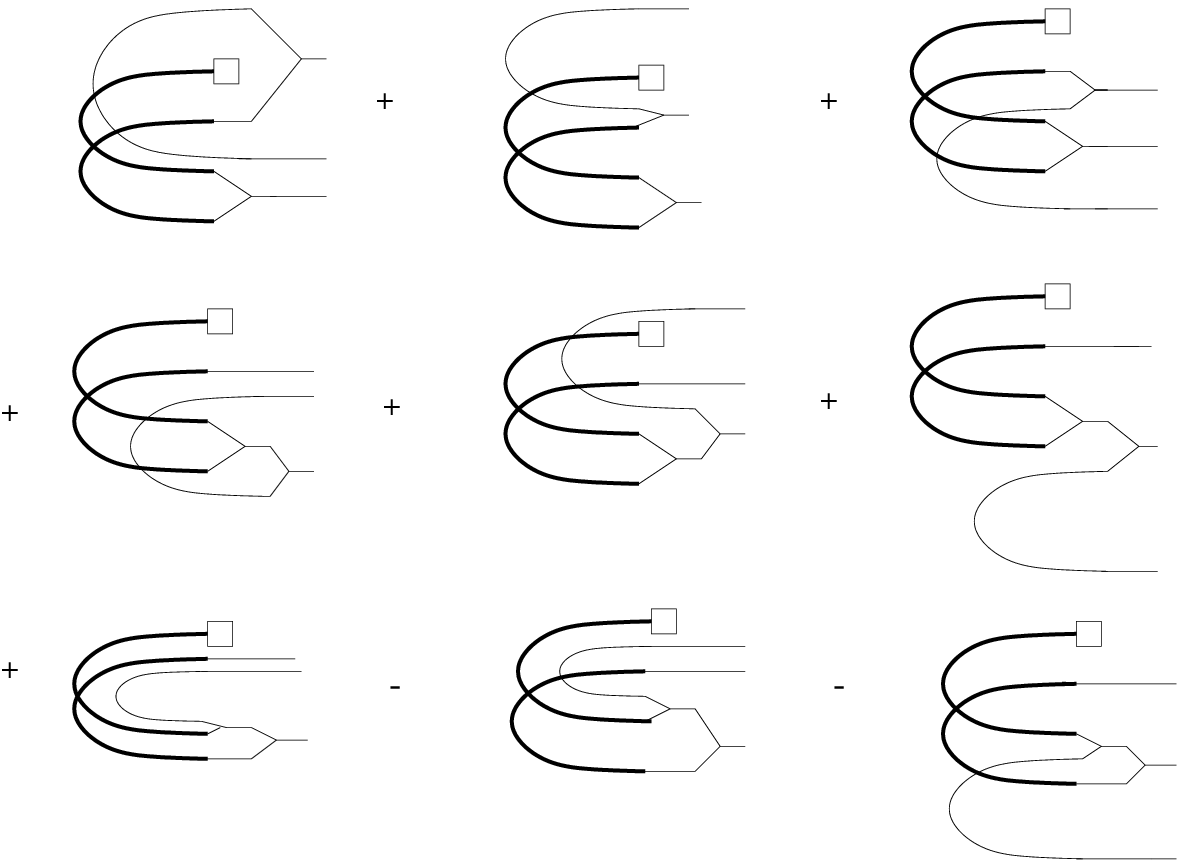}}
$$
The sums of terms in each column is zero by Jacobi Identity. This concludes
the proof of Theorem~\ref{T:g01}.
\begin{flushright}
$\square\;\;\;\;\;\;\;\;$
\end{flushright}
\section{Classification of dynamical r-matrices}
\paragraph{}Let $\g$ be a simple algebra. In that case,~(\ref{E:08}) becomes
\begin{equation}\label{E:09}
r(\lambda) + r^{21}(\lambda)=\epsilon \Omega.
\end{equation}
\hbox to1em{\hfill}We will classify all solutions of
equations~(\ref{E:07},\ref{E:06},\ref{E:09}) when $\epsilon \neq 0$ and
when $\l$ contains a semisimple regular element. In particular, in this
case, the centralizer $\h$ of $\l$ is the unique Cartan subalgebra
containing $\l$. Notice that we can assume that $\epsilon =1$ ( since the
assignement $r(\lambda) \to \epsilon r(\epsilon \lambda)$ is a gauge
transformation of~(\ref{E:06})). We can also assume that the restriction of
$(\,,\,)$ to $\l$ is nondegenerate. Indeed, for any dynamical r-matrix, we
can replace $\l$ by the largest subspace of $\h$ for which $r$ is
invariant, and such a subspace is real. Let $\h_{0}$ be the orthogonal
complement of $\l$ in $\h$ and let $i: \l \hookrightarrow \h$ be the
inclusion map. We will also write $(\,,\,)$ for the induced bracket on $\l^*$. Let $\Omega_{\h_0}$ denote the Casimir element of the
restriction of $(\,,\,)$ to $\h_0$.
\subsection{Statement of the theorem}
\paragraph{}Let $\g=\n_{+}\oplus \h \oplus \n_{-}$ be a polarization of $\g$.\\
{\bf{Definition:}} A generalized Belavin-Drinfeld triple is a triple
$(\Gamma_{1},\Gamma_{2},\tau)$ where $\Gamma_{1},\Gamma_{2} \subset \Pi$,
and $\tau:\Gamma_{1}\stackrel{\sim}{\to}\Gamma_{2}$ is a norm-preserving
bijection.\\
\hbox to1em{\hfill}In other terms, in a generalized Belavin-Drinfeld
triple, we drop the nilpotency condition. We will say that a generalized
Belavin-Drinfeld triple is $\l$-graded if $\tau$ preserves the
decomposition of $\g$ in $\l$-weight spaces. If
$(\Gamma_{1},\Gamma_{2},\tau)$ is a generalized Belavin-Drinfeld triple, we
will denote by $\Gamma_{3}$ the largest subset of $\Gamma_{1} \cap
\Gamma_{2}$ which is stable under $\tau$, and
$\tilde{\Gamma}_{1}=\Gamma_{1} \backslash \Gamma_{3}$,
$\tilde{\Gamma}_{2}=\Gamma_{2} \backslash \Gamma_{3}$. It is clear that
$(\tilde{\Gamma}_{1},\tilde{\Gamma}_{2},\tau)$ is a Belavin-Drinfeld
triple. As before, for each choice of Chevalley generators
$(e_{\alpha},f_{\alpha},h_{\alpha})_{\alpha \in \Gamma_{i}}$, the map
$\tau$ induces isomorphisms $\g_{\tilde{\Gamma}_{1}} \to
\g_{\tilde{\Gamma}_{2}}$ and $\tau: \g_{\Gamma_{3}} \to \g_{\Gamma_{3}}$.\\
\hbox to1em{\hfill}For $\lambda \in \l^*$, consider the map:
\begin{align*}
K(\lambda):\n_{+}({\Gamma_{1}}) &\to \n_{+}({\Gamma_{2}})\\
	e_{\alpha} &\mapsto \e e_{\alpha} +  e^{-(\alpha,\lambda)}
\frac{\tau}{1-e^{-(\alpha,\lambda)}\tau}(e_{\alpha}).
\end{align*}
Notice that we have
$$K(\lambda)(e_{\alpha})=\frac{1}{2}e_{\alpha} + \sum_{n >0}
e^{-n(\alpha,\lambda)}\tau^{n}(e_{\alpha}).$$
This sum is finite for $\alpha \not\in \left< \Gamma_3 \right>$.
\begin{theo}Let $\g$ be a simple Lie algebra with nondegenerate invariant bilinear form $(\,,\,)$, $\l \subset \g$ a commutative
subalgebra containing a regular semisimple element on which $(\,,\,)$ is nondegenerate, $\h$ the Cartan
subalgebra containing $\l$ and $\h_{0}$ the orthogonal complement of $\l$
in $\h$. Then
\begin{enumerate}
\item Any dynamical r-matrix is gauge-equivalent to a dynamical r-matrix
$\tilde{r}$ such that
\begin{equation}\label{E:09.25}
 \tilde{r}(\lambda)-\tilde{r}(\lambda)^{21} \in (\l^{\perp})^{\otimes
2}=(\bigoplus_{\alpha \neq 0}\g_{\alpha} \oplus \h_{0})^{\otimes 2}.
\end{equation}
\item Let $(\Gamma_{1},\Gamma_{2},\tau)$ be an $\l$-graded generalized
Belavin-Drinfeld triple and let $(e_{\alpha}, f_{\alpha},
h_{\alpha})_{\Gamma_{i}}$ be a choice of Chevalley generators. Let
$r_{\h_{0},\h_{0}} \in \h_{0} \otimes \h_{0}$ satisfy the equation
\begin{equation}\label{E:09.5}
(\tau(\alpha) \otimes 1)r_{\h_0,\h_0} + (1 \otimes
\alpha)r_{\h_0,\h_0}=\frac{1}{2}((\alpha + \tau(\alpha)) \otimes
1)\Omega_{\h_0}.
\end{equation}
Then
$$r(\lambda)=\frac{1}{2}\Omega + r_{\h_0,\h_0} + \sum_{\alpha \in
\left<\Gamma_{1}\right> \cap \Delta_+} K(\lambda)(e_{\alpha}) \wedge
e_{-\alpha}+ \sum_{\alpha \in \Delta_+,\, \alpha \not\in \left< \Gamma_1
\right>} \frac{1}{2}e_{\alpha} \wedge e_{-\alpha}$$
is a solution the CDYBE satisfying~(\ref{E:09.25}).
\item Any solution of the CDYBE  satisfying~(\ref{E:09.25}) is of the above
type for a suitable polarization of $\g$.
\end{enumerate}
\end{theo}
The proof of this theorem will occupy the rest of this section. Our methods
are greatly inspired by the paper \cite{BD}. Notice that 1. follows from
Theorem~\ref{T:g01}, but we will describe the gauge transformations
explicitely in this case.
\paragraph{Notations:}Let $\Delta \subset \h^*$ be the root system of $\g$
with respect to $\h$ and set $\Delta_{\l}=i^*(\Delta) \subset \l^*$. We
will denote by $\g_{\oa}$ the weight subspace associated to $\oa=i^*(\alpha) \in
\Delta_{\l}$, and we set $\g_{\overline{0}}=\h_{0}$. It is clear that
$$\g_{\overline{\alpha}}=\bigoplus_{\beta \in
\Delta,\;i^*(\beta)=\overline{\alpha}} \g_{\beta}$$
 In particular, $(\,,\,)$ is a pairing $\g_{\overline{\alpha}} \times
\g_{-\overline{\alpha}} \to \C$.\\
\hbox to1em{\hfill}A vector space $V \subset \g$ will be called $\h$-graded
(resp. $\l$-graded) if it is an $\h$-submodule (resp. $\l$-submodule) of
$\g$. Finally, let $\Omega'\in (\l^{\perp})^{\otimes 2}$ denote the Casimir
(inverse element) of the restriction of $(\,,\,)$ to $\l^{\perp}=\h_{0}
\bigoplus \g_{\overline{\alpha}}$.
\paragraph{}Now let $r:\l^*\supset D \to (\g \otimes \g)^{\l}$ be a formal power series satisfying~(\ref{E:09}) (with $\epsilon =1$).
By~(\ref{E:07}), we can write
\begin{equation}\label{E:10}
r(\lambda)=\frac{1}{2}\Omega+r_{\l,\l}(\lambda)+r_{\l,\h_{0}}(\lambda) +
r_{\h_{0},\l}(\lambda) + (\varphi(\lambda) \otimes 1)\Omega',
\end{equation}
where $r_{\l,\l}(\lambda) \in \l \otimes \l$, $r_{\l,\h_0}(\lambda) \in \l
\otimes \h_0$, $r_{\h_0,\l}(\lambda) \in \h_0 \otimes \l$ and where
$\varphi(\lambda)\in \mathrm{End}\,\big(\h_{0} \bigoplus
\g_{\overline{\alpha}}\big)$ is a sum of maps
$\varphi_{\overline{\alpha}}(\lambda)\in
\mathrm{End}\,(\g_{\overline{\alpha}})$. By the unitarity condition,
$r_{\l,\l}(\lambda) \in \Lambda^2 \l$,
$r_{\l,\h_{0}}(\lambda)=-r_{\h_{0},\l}^{21}(\lambda)$ and
$\varphi_{-\overline{\alpha}}(\lambda)=-\varphi_{\overline{\alpha}}^{*}(\lambda)
$.\\
\hbox to1em{\hfill}With these notations, the CDYBE splits into $4$
components: the $\l \otimes \l \otimes \l$-part, the $\l \otimes \l \otimes
\h_{0}$-part, the $\l\otimes \g_{\oa} \otimes \g_{-\oa}$-part and the
$\g_{\oa} \otimes \g_{\ob} \otimes \g_{\oj}$-part where $\oa + \ob + \oj=0$.
\begin{itemize}
\item The $\l \otimes \l \otimes \l$-part: let us set
$r_{\l,\l}=\sum_{i,j}C_{i,j}(\lambda) x_{i} \otimes x_{j}$. This part of
the CDYBE can then be written:
\begin{equation}\label{E:11}
\frac{\partial C_{j,k}}{\partial x_{i}}+\frac{\partial C_{k,i}}{\partial
x_{j}} + \frac{\partial C_{i,j}}{\partial x_{k}}=0  \qquad \forall\;i,j,k
\end{equation}
and says that $\sum_{i,j}C_{i,j}dx_{i} \wedge dx_{j}$ is a closed 2-form.
\item The $\l \otimes \l \otimes \h_{0}$-part: let us set
$r_{\l,\h_{0}}=\sum_{i,j} D_{i,j}(\lambda) x_{i} \otimes y_{j}$ for some
basis $(y_{j})$ of $\h_{0}$. This part of the CDYBE is
\begin{equation}\label{E:12}
\frac{\partial D_{i,j}}{\partial x_{k}}=\frac{\partial D_{k,j}}{\partial
x_{i}} \qquad \forall\;i,k,j
\end{equation}
and says that for any $j$, $\sum_{i}D_{i,j}(\lambda)dx_{i}$ is a closed 1-form.
\end{itemize}
\hbox to1em{\hfill}Since $r$ is defined on a polydisc, the above forms are
exact. Let $f:\,D \to \h_0$ be such that $df(\lambda)=\sum_i
D_{i,j}(\lambda)dx_i \otimes y_j$ and let $\xi$ be a 1-form on $D$ such that $d\xi = \sum_{i,j}C_{i,j}dx_i \wedge dx_j$. Then $\xi$ defines a function $\overline{\xi}: D \to \l$. The gauge transformation which should be applied to $r$
to make it satisfy~(\ref{E:09.25}) is easily seen to be the following:
$$r(\lambda) \mapsto r(\lambda)^{g}=\frac{1}{2}\Omega +
(e^{-ad\;f(\lambda)}\varphi(\lambda)e^{ad\;f(\lambda)} \otimes 1) \Omega'$$
where $g(\lambda)=e^{f(\lambda)}e^{-\overline{\xi}(\lambda)}$.
\paragraph{}Thus, we can assume that $r_{\l,\l}=r_{\l,\h_0}=0$, in which case the remaining components of the CDYBE can be written in the following
way:
\begin{itemize}
\item The $\l\otimes \g_{\oa} \otimes \g_{-\oa}$-part:
\begin{equation}\label{E:13}
d\varphi_{\oa} + (\varphi^2_{\oa}-\frac{1}{4})dh_{\oa}=0.
\end{equation}
In particular, $r_{\h_0,\h_0} \in \Lambda^2 \h_0$ is constant.
\item The $\g_{\oa} \otimes \g_{\ob} \otimes \g_{\oj}$-part where $\oa +
\ob + \oj=0$:
\begin{equation}\label{E:14}
\Lambda \big( \varphi_{\oa} \otimes \varphi_{\ob} \otimes 1 + \varphi_{\oa}
\otimes 1 \otimes \varphi_{\oj} + 1 \otimes \varphi_{\ob} \otimes
\varphi_{\oj}  + \frac{1}{4}Id \big)=0
\end{equation}
where $\Lambda:\;\g_{\oa} \otimes \g_{\ob} \otimes \g_{\oj} \to \C$, $x
\otimes y \otimes z \mapsto ([x,y],z)$.
\end{itemize}
\hbox to1em{\hfill}This set of equations is sufficient by skew-symmetry of
the CDYBE.
\subsection{The Cayley transform}
\paragraph{}Let us set $A_{\pm}=\mathrm{Im}\,(\ph (\lambda) \pm \frac{1}{2})$,
$I_{\pm}=\mathrm{Ker}\,(\ph (\lambda) \mp \e)$. Notice that, by (\ref{E:13}), $A_{\pm}$ and $I_{\pm}$ are indeed independent of $\lambda$. Furthermore, $A_{\pm},\,I_{\pm}$ are
$\l$-graded by the weight-zero condition, $I_{\pm} \subset A_{\pm}$ and
$A_{\pm}=I_{\pm}^{\perp}$ by the unitarity condition. Notice also that
$A_{+} + A_{-} \oplus \l=\g$. Now consider
$$
\psi(\lambda)=\frac{\ph-\e}{\ph + \e}:\;A_{+}/I_{+} \to A_{-}/I_{-}.
$$
Extend $\psi(\lambda)$ to ${\psi}(\lambda): \l \oplus A_{+}/I_{+} \to \l \oplus A_{-}/I_{-}$ by setting ${\psi}_{|\l}=Id$.
It is clear that ${\psi}$ is a well-defined linear isomorphism. The following proposition is crucial:
\begin{prop} The maps $\ph_{\oa}$ satisfy~(\ref{E:13},\ref{E:14}) if and
only if the following hold:
\begin{enumerate}
\item[(i)] $A_{\pm}\oplus \l$ is a subalgebra of $\g$ and $I_{\pm} \oplus
\l$ is an ideal of $A_{\pm} \oplus \l$.
\item [(ii)] there exists a (constant) map ${\psi}_{0}:\;\l \oplus A_{+}/I_{+} \to \l \oplus
A_{-}/I_{-}$ such that
$$\psi(\lambda)_{|\g_{\oa}}=e^{-(\oa,\lambda)}{\psi}_{0|\g_{\oa}}.$$
\item[(iii)] The map ${\psi}_0$ is a Lie algebra map:
\begin{equation}\label{E:17}
[\psi_0(x),\psi_0(y)]=\psi_0[x,y].
\end{equation}
\end{enumerate}
\end{prop}
{\bf{Proof:}} Assume that $\ph$ satisfies~(\ref{E:13},\ref{E:14}) and let
$a \in \g_{\oa},\,b \in \g_{\ob},\,c \in \g_{\oj}$ with $\oa + \ob +\oj=0$.
From~(\ref{E:14}), we have
\begin{align*}
\big( [(\ph_{\oa}+ \e)a,(\ph_{\ob} + \e)b],c\big)+& \big([a,(\ph_{\ob} +
\e)b],(\ph_{\oj} - \e)c\big)\\
 +& \big([(\ph_{\oa}- \e)a,b],(\ph_{\oj} - \e)c\big)=0.
\end{align*}
Since $\ph_{\oj}=-\ph_{-\oj}^{*}$, and $(\,,\,)$ is a nondegenerate pairing
$\g_{\oj} \otimes \g_{-\oj} \to \C$, this implies that $A_{+} \oplus \l$ is
a Lie subalgebra of $\g$. Note that the term in $\l$ is necessary here
since $[\g_{\oa},\g_{-\oa}] \not\subset \g_{\overline{0}}=\h_{0}$, but is
not consequential as $A_{+}$ is $\l$-graded. The second claim of (i)
follows from the relation
\begin{align*}
\big( [(\ph_{\oa}- \e)a,(\ph_{\ob} -\e)b],c\big)&+ \big([a,(\ph_{\ob} +
\e)b],(\ph_{\oj} + \e)c\big)\\
 &+ \big([(\ph_{\oa}- \e)a,b],(\ph_{\oj} + \e)c\big)=0.
\end{align*}
The proof is the same for $A_{-}$ and $I_{-}$. The equivalence of (ii)
and~(\ref{E:13}) follows from the equality
\begin{equation*}
\begin{split}
d{\psi}_{|\g_{\oa}}&= \frac{d\ph_{\alpha} (\ph_{\alpha} + \frac{1}{2})-(\ph_{\alpha}-\frac{1}{2})d\ph_{\alpha}}{(\ph_{\alpha}+ \frac{1}{2})^2}\\
&=- \frac{(\ph_{\alpha}^2-\frac{1}{4})}{(\ph_{\alpha} + \frac{1}{2})^2}dh_{\oa}\\
&=-(\oa, \lambda){\psi}_{|\g_{\oa}}.
\end{split}
\end{equation*}
where we used (\ref{E:13}). Finally it follows
from~(\ref{E:14}) that
$$(\ph_{\overline{\alpha + \beta}} - \e) \big([(\ph_{\oa}+ \e)a,(\ph_{\ob}
+ \e)b]\big)=(\ph_{\overline{\alpha + \beta}} + \e) \big([(\ph_{\oa}-
\e)a,(\ph_{\ob} - \e)b]\big).$$
This implies (iii).\\
\hbox to1em{\hfill}Conversely, if (i-iii) are satisfied then for any $x \in
\g_{\oa},\,y\in \g_{\ob}$ ($\oa+\ob \neq 0$) there exist $ z\in
\g_{\overline{\alpha+\beta}}$ such that
$$[(\ph_{\oa}-\e)x,(\ph_{\ob}-\e)y]=(\ph_{\overline{\alpha+\beta}}-\e)z.$$
Since $\psi$ is a Lie algebra map,
$[(\ph_{\oa}+\e)x,(\ph_{\ob}+\e)y]-(\ph_{\overline{\alpha+\beta}}+\e)z \in
\mathrm{Ker}\,(\ph_{\overline{\alpha + \beta}} -\e)$. Subtracting, we obtain
$$[(\ph_{\oa}+\e)x,y]+ [x,(\ph_{\ob}+\e)y] -[x,y] -z \in \mathrm{Ker}\,(\ph_{\overline{\alpha + \beta}}
-\e).$$
Applying $(\ph -\e)$ and dropping the indices, we have
$$(\ph - \frac{1}{2})\Big([(\ph + \frac{1}{2})x,y] + [x,(\ph + \frac{1}{2})y] -[x,y] \Big)=[(\ph -\frac{1}{2})x,(\ph - \frac{1}{2})y].$$
Thus,
$$[(\ph + \frac{1}{2})x,(\ph + \frac{1}{2})y]-(\ph + \frac{1}{2}) \Big( [(\ph - \frac{1}{2})x,y] + [x,(\ph + \frac{1}{2})y] \Big)=0.$$
which is equivalent to (\ref{E:14}). 
\begin{flushright}
$\square\;\;\;\;\;\;\;\;$
\end{flushright}
\paragraph{}We will call the triple $(A_{+},A_{-},\psi_{0})$ the Cayley
transform of $\ph$. We are now reduced to the classification of all triples
satisfying (i-iii) and which arise as a Cayley transform (Cayley triples).
\subsection{Classification of Cayley triples}
\paragraph{}Let $(A_{+},A_{-},\psi_{0})$ be a Cayley triple. If $\g=\n_{+}
\oplus \h \oplus \n_{-}$ is a polarization of $\g$ and $\Gamma \subset
\Pi(\n_+)$ we will denote by $\q^{+}_{\Gamma}$ (resp. $\q_{\Gamma}^{-}$)
the subalgebra generated by $\n_{+}$ and $\g_{-\alpha},\,\alpha \in \Gamma$
(resp. generated by $\n_{-}$ and $\g_{\alpha},\,\alpha \in \Gamma$). We
denote by $\p^{\pm}_{\Gamma}=\h + \q_{\Gamma}^{\pm}$ the parabolic
subalgebras associated to $\Gamma$.
\begin{prop}\label{P:01}There exists a polarization $\g=\n_{+}^1 \oplus \h
\oplus \n_{-}^1$, two subsets $\Gamma_{+},\Gamma_- \subset \Pi(\n_{+}^1)$
and two vector spaces $V_{+},V_{-} \subset \h$ with $V_{\pm}^{\perp}
\subset V_{\pm}$ such that
$$\l \oplus A_{+}=\q_{\Gamma_{+}}^{+} \oplus V_{+},\qquad \l \oplus
A_{-}=\q_{\Gamma_-}^{-} \oplus V_{-}$$
\end{prop}
{\bf{Proof:}} Notice that $(\l \oplus A_{+})^{\perp}=I_{+} \subset \l
\oplus A_{+}$. It is known, (c.f [Bou, chap.VIII,\S 10, Thm. 1] or
\cite{BD}), that this implies that $\l \oplus A_{+}=\tilde{\q}_{\Gamma}^{+}
\oplus \tilde{V}_{+}$ for {\it{some}} polarization $\g=\n_{+}' \oplus \h'
\oplus \n_{-}'$. Similarly, $\l \oplus A_{-}=\tilde{\q}_{\Gamma'}^{-}
\oplus \tilde{V}_{-}$ for some polarization $\g=\n_{+}'' \oplus \h''\oplus
\n_{-}''$. Moreover, $\l$ acts semisimply on $A_{\pm}$ so $\l \subset \h'$, $\l \subset \h''$.
But $\l$ contains a regular element, thus $\l=\h'=\h''$.
Proposition~\ref{P:01} is now an easy consequence of the following lemma:
\begin{lem}\label{L:01} Let $\g$ be a simple Lie algebra and $\h$ a Cartan
subalgebra. Let $\a_{1}$ and $\a_{2}$ be two parabolic subalgebras
containing $\h$ such that $\a_{1} + \a_{2}=\g$. Then there exists a
polarization $\g=\n_{+} \oplus \h \oplus \n_{-}$ and $\Gamma_{+},\Gamma_-
\subset \Pi$ such that $\a_{1}=\p_{\Gamma_{+}}^{+}$ and
$\a_{2}=\p_{\Gamma_-}^{-}$.
\end{lem}
{\bf{Proof:}} Let $\n_{+} \oplus \h \oplus \n_{-}$ be a polarization of
$\g$ such that $\b_{+} \subset \a_{1}$ and for which $\mathrm{dim}\,
(\n_{+} \cap \a_{2})$ is minimal. We claim that $\b_{-} \subset \a_{2}$.
Suppose on the contrary that there exists a simple root $\alpha \in \Pi$
such that $\g_{-\alpha} \not \subset \a_{2}$. Then $\g_{-\alpha} \subset
\a_{1}$ since $\a_{1} + \a_{2}=\g$ and $\g_{\alpha} \subset \a_{2}$ since
$\a_{2}$ is parabolic. But then $s_{\alpha}\n_{+} \oplus \h \oplus
s_{\alpha}\n_{-}$ is a polarization of $\g$ for which $s_{\alpha}\b_{+}
\subset \a_{1}$ and $\mathrm{dim}\,(s_{\alpha}\n_{+} \cap \a_{2}) <
\mathrm{dim}\,(\n_{+} \cap \a_{2})$. Contradiction.
\begin{flushright}
$\square\;\;\;\;\;\;\;\;$
\end{flushright}
\hbox to1em{\hfill}In particular, $A_{\pm}$, $I_{\pm}$ are all $\h$-graded and
\begin{align*}
I_{+}&=(\q_{\Gamma_{+}}^+ \oplus V_{+})^{\perp}= \bigoplus_{\alpha \in
\Delta_{+} \backslash \left<\Gamma_{+}\right>}\g_{\alpha} \oplus
(V_{+}^{\perp} \cap \h_{0}),\\
I_{-}&=(\q_{\Gamma_-}^- \oplus V_{-})^{\perp}=\bigoplus_{\alpha \in
\Delta_{-} \backslash \left<\Gamma_-\right>}\g_{\alpha}
\oplus(V_{-}^{\perp} \cap \h_{0}).
\end{align*}
Thus $A_{+} / I_{+}= \g_{\Gamma_{+}} \oplus V_{1}$ and
$A_{-}/I_{-}=\g_{\Gamma_-} \oplus V_{2}$ for some suitable $V_{1},V_{2}
\subset \h_{0}$.
\paragraph{}Let $L_{\pm \e}(\lambda)$ be the generalized eigenspace of
$\ph(\lambda)$ associated to $\pm \e$. Since $\ph$ is a solution of an
ordinary differential equation with stationary points at $\e,-\e$, $L_{\pm
\e}(\lambda)$ is independent of $\lambda$ and we will simply denote it by
$L_{\pm \e}$. Similarly, let $L'$ be the sum of all other generalized
eigenspaces so that $\g=\l \oplus L_{\e} \oplus L' \oplus L_{-\e}$.
\begin{prop}\label{P:02}There exists a polarization $\g=\overline{\n}_{+}
\oplus \h \oplus \overline{\n}_{-}$ and a subset $\Gamma_{3} \subset
\Pi(\overline{\n}_{+})$ such that $L_{\pm \e} \subset \overline{\b}_{\pm}$,
$L' \subset \g_{\Gamma_{3}} + \h$ and $\ph(\overline{\n}_{+}) \subset
\overline{\n}_+$.
\end{prop}
{\bf{Proof:}} We will construct a polarization satisfying the above
conditions in several steps.
\begin{lem}\label{L:02} We have:
\begin{enumerate}
\item[(i)] $\l \oplus L_{\pm \e}$ is an $\h$-graded solvable subalgebra,
\item[(ii)] $\l \oplus L'$ is an $\h$-graded subalgebra,
\item[(iii)] we have $[L_{\pm \e},L'] \subset \l \oplus L_{\pm \e}$.
\end{enumerate}
\end{lem}
{\bf{Proof:}} this follows from the proofs of Lemma 12.3 and Theorem 12.6
in \cite{BD}.
\paragraph{}Notice that $L_{\pm \e} \not\subset \b_{\pm}^1$ in general. We
first construct a polarization $\g=\n_+^2 \oplus \h \oplus \n_-^2$ such that $L_{\pm \e} \subset \b_{\pm}^2$. We have
$I_{\pm} \subset L_{\pm \e}$. Notice that $L_{\e} \cap \n_{-}^1 \subset
\g_{\Gamma_{+}} \cap \g_{\Gamma_-}=\g_{\Gamma_{+} \cap \Gamma_-}$ since
$\n_{-}^1 \subset (\g_{\Gamma_-} \oplus I_{-})$ and $L_{\e}$ is solvable. Similarly, $L_{-\e} \cap
\n_{+}^1 \subset \g_{\Gamma_{+} \cap \Gamma_-}$. Moreover, by
Lemma~\ref{L:02}, $\l \oplus (L_{\e} \cap \g_{\Gamma_{+} \cap \Gamma_-})$
and $\l \oplus (L_{-\e} \cap \g_{\Gamma_{+} \cap \Gamma_-})$ are disjoint,
solvable, $\h$-graded subalgebras. By lemma~\ref{L:01} it follows that
there exists an element $s$ of the group $W_{\Gamma_{+} \cap \Gamma_-}$
such that $$\l \oplus (L_{\pm \e} \cap \g_{\Gamma_{+} \cap \Gamma_-})
\subset s\b_{\pm}^1.$$
 Notice that $s$ permutes elements of $\Delta^+ \backslash
\left<\Gamma_{+} \cap \Gamma_-\right>$, leaving it globally unchanged.
Thus, $\l \oplus L_{\pm \e} \subset s\b_{\pm}^1$. Set $\n_{\pm}^2=s
\n_{\pm}^1$.\\
\hbox to1em{\hfill}Now we construct a polarization of $\g$ satisfying the
other conditions of proposition~\ref{P:02}. Recall that $\l \oplus L
\subset \g_{\Gamma_{+} \cap \Gamma_-} + (V_{1} \cap V_{2})$. Thus
$$\big(L' \cap \n_{+}^2\big) \oplus \big(L_{\e} \cap \n_{+}^2(\Gamma_{+}
\cap \Gamma_-)\big)=\n_{+}^2(\Gamma_{+} \cap \Gamma_-).$$
\hbox to1em{\hfill}Since $[L',L_{\e}]\subset \l \oplus L_{\e}$ by
Lemma~\ref{L:02},(iii), $L_{\e} \cap \n_{+}^2(\Gamma_{+} \cap \Gamma_-)$ is
an ideal of $\n_{+}^2(\Gamma_{+} \cap \Gamma_-)$. But $L' \cap \n_{+}^2$ is
a subalgebra. It is easy to see that this implies that $L' \cap \n_{+}^2$
is generated by a set of simple root subspaces of $\n_{+}^2(\Gamma_{+} \cap
\Gamma_-)$, i.e $L' \cap \n_{+}^2=\n_{+}^2(\Gamma)$ for some $\Gamma
\subset \Pi(\n_{+}^2)$. Moreover, the restriction of $(\,,\,)$ to $L'$ is
nondegenerate, hence $L' \cap \n_{-}^2=\n_{-}^2(-\Gamma)$. Thus
$$\l \oplus \g_{\Gamma} \subset \l \oplus L' \subset \l\oplus \g_{\Gamma} +
(V_{1} \cap V_{2}).$$
\hbox to1em{\hfill}Since $\ph(\lambda) + \e$ is invertible in $L'$,
$\psi(\lambda)$ is a well-defined operator $L' \to L'$,
satisfying~(\ref{E:17}), and $\psi(\lambda) (\h_{0} \cap L')
\subset \h_{0} \cap L'$. Now, $\l$ contains a regular element. Thus there
exists a polarization of $\g$ compatible with the $\l$-weight
decomposition. This induces a polarization of $\g_{\Gamma}$, compatible
with the $\l$-weight decomposition of $\g_{\Gamma}$. Hence, there exists
$s' \in W_{\Gamma} \subset W$ such that $\psi_{0|\g_{\Gamma}}$ is
compatible with the polarization $s'\n_{+}^2 \oplus \h \oplus s'\n_{-}^2$.
Since $s'$ leaves $\Delta_+ \backslash \left<\Gamma\right>$ globally
unchanged, the polarization $\g=\overline{\n}_{+} \oplus \h \oplus
\overline{\n}_{-}$ with $\overline{\n}_{\pm}=s'\n_{\pm}^2$ and
$\Gamma_3=s'\Gamma$ satisfies the requirements of proposition~\ref{P:02}.
\begin{flushright}
$\square\;\;\;\;\;\;\;\;$
\end{flushright}
\paragraph{}To sum up, we have shown that there exists a polarization
$\g=\overline{\n}_{+} \oplus \h \oplus \overline{\n}_{-}$, compatible with
$\ph$, subsets $\Gamma_1=s's\Gamma_+$, $\Gamma_2=s's\Gamma_-$ and
$\Gamma_3\,\subset \Pi(\overline{\n}_{+})$ such that $(A_{+}/I_{+})\cap
\n_{+}=\overline{\n}_{+}({\Gamma_{1}})$, $A_{-}\cap
\n_{+}=\overline{\n}_{+}({\Gamma_{2}})$ and $L' \cap
\n_{+}=\overline{\n}_{+}(\overline{\Gamma_{3}})$.\\
\hbox to1em{\hfill}The map $\psi_{0}$ now restricts to a Lie algebra
isomorphism $\overline{\n}_{+}(\Gamma_{1}) \to
\overline{\n}_{+}(\Gamma_{2})$. This isomorphism maps weight spaces to
weight spaces as $\psi_0$ preserves $\h_0$ and $\ph$ is $\l$-invariant.
Define $\tau: \Gamma_{1} \to \Gamma_{2}$ by
$\psi_{0}(\g_{\alpha})=\g_{\tau(\alpha)}$. It is a norm-preserving
bijection. Thus $(\Gamma_1,\Gamma_2,\Gamma_3)$ is a generalized
Belavin-Drinfeld triple. It is clear that $\Gamma_{3}$ is the largest
subset of $\Gamma_{1} \cap \Gamma_{2}$ stable under $\tau$, and that
$\psi_0: \overline{\n}_+(\Gamma_3) \to \overline{\n}_+(\Gamma_3)$ is a Lie algebra isomorphism.
Finally, it is easy to see that the map $\varphi$ is obtained from this
data by formulas
\begin{alignat*}{2}
\ph(\lambda)(e_{\alpha})&=\frac{1}{2}e_{\alpha}& \qquad &(\alpha \not \in
\left< \Gamma_1 \right>)\\
\ph(\lambda)(e_{\alpha})&=\frac{1}{2}e_{\alpha} +
\frac{\psi_0}{1-e^{(\alpha,\lambda)}\psi_0}(e_{\alpha})& \qquad &(\alpha
\in \left< \Gamma_1 \right>)
\end{alignat*}
\paragraph{}Conversely, it is clear how to construct from a generalized
Belavin-Drinfeld triple $(\Gamma_1,\Gamma_2,\tau)$ the subalgebras
$\n_+(\Gamma_1)$, $\n_+(\Gamma_2)$, $\n_+(\Gamma_3)$ and, for each choice
of Chevalley generators, a Lie algebra isomorphism $\psi_0: \n_+(\Gamma_1)
\to \n_+(\Gamma_2)$, and the map $\ph(\lambda)$. Condition~(\ref{E:09.5})
on the $\h_0 \otimes \h_0$-part comes from~(\ref{E:14})-see \cite{BD}.

\section{Examples}
\subsection{Constant r-matrices}
\paragraph{}Our results imply the following:
\begin{cor} A dynamical r-matrix associated to a generalized
Belavin-Drinfeld triple $(\Gamma_1,\Gamma_2,\tau)$ is gauge equivalent to a
constant r-matrix if and only if $\Gamma_3=\emptyset$.
\end{cor}
\subsection{$\h$-invariant dynamical r-matrices}
\paragraph{}When $\l=\h$, our classification coincides with that given in
\cite{EV}: the only $\h$-graded generalized Belavin-Drinfeld triple is of
the form $(\Gamma,\Gamma, \tau=Id)$. The dynamical r-matrices obtained are
then (up to gauge transformations and choice of Chevalley generators):
$$r(\lambda) = \frac{\Omega}{2} + \sum_{\alpha \in \Delta_+,\,\alpha
\not\in \left<\Gamma\right>} \frac{1}{2}e_{\alpha} \wedge e_{-\alpha} +
\sum_{\alpha \in \left<\Gamma\right> \cap\Delta_+}
\frac{1}{2}coth((\frac{1}{2}(\alpha, \lambda)) e_{\alpha} \wedge
e_{-\alpha}.$$
\subsection{Example for $\mathfrak{sl}_{3}$ and $\mathfrak{sl}_n$}
\paragraph{}The first nontrivial example is for $\g=\mathfrak{sl}_3$: fix a
polarization $\g=\h \oplus \bigoplus_{\gamma \in \Delta} \g_{\gamma}$ where
$\Delta^+=\{\alpha,\beta,\alpha + \beta\}$ and set $\l=\C h_{\rho}$. Consider
the generalized Belavin-Drinfeld triple with
$\Gamma_{1}=\Gamma_2=\{\alpha,\beta\}$ and $\tau: \alpha \mapsto
\beta,\,\beta \mapsto \alpha$. In this case, we can choose the map $\psi_0$
to be the following
\begin{alignat*}{3}
e_{\alpha} &\mapsto e_{\beta}, & \quad h_{\alpha} & \mapsto h_{\beta}, &
\quad e_{-\alpha}&\mapsto e_{-\beta}\\
e_{\beta} &\mapsto e_{\alpha}, & \quad h_{\beta} & \mapsto h_{\alpha}, &
\quad e_{-\beta}&\mapsto e_{-\alpha}\\
e_{\alpha + \beta} &\mapsto -e_{\alpha + \beta}, & \quad &  & \quad
e_{-\alpha-\beta}&\mapsto -e_{-\alpha-\beta}.
\end{alignat*}

The corresponding dynamical r-matrix is given by:
\begin{equation}
\begin{split}
r(\lambda)= \frac{\Omega}{2} &+ r_{\h_0,\h_0} + \frac{1}{2}coth(\alpha,
\lambda)e_{\alpha} \wedge e_{-\alpha} + \frac{1}{2}coth(\beta,
\lambda)e_{\beta} \wedge e_{-\beta}\\
 &+ \frac{1}{2}th(\alpha + \beta, \lambda)e_{\alpha+\beta} \wedge
e_{-\alpha-\beta} + \frac{1}{2 sinh(\alpha, \lambda)}e_{\beta} \wedge
e_{-\alpha}\\
& + \frac{1}{2 sinh((\alpha, \lambda))}e_{\alpha} \wedge e_{-\beta}.
\end{split}
\end{equation}
This dynamical r-matrix is gauge-equivalent to the dynamical r-matrix
\begin{equation}
\begin{split}
\tilde{r}(\lambda)= \frac{\Omega}{2} &+ r_{\h_0,\h_0} +r_{\l,\h_0} -
r_{\l,\h_0}^{21} + \frac{1}{2}coth(\alpha, \lambda)e_{\alpha} \wedge
e_{-\alpha} + \frac{1}{2}coth(\beta, \lambda)e_{\beta} \wedge e_{-\beta}\\
 &+ \frac{1}{2}th(\alpha + \beta, \lambda)e_{\alpha+\beta} \wedge
e_{-\alpha-\beta} + \frac{e^{(\alpha,\lambda)}}{2 sinh(\alpha,
\lambda)}e_{\beta} \wedge e_{-\alpha} \\
&+ \frac{e^{-(\alpha,\lambda)}}{2 sinh(\alpha, \lambda)}e_{\alpha} \wedge
e_{-\beta}.
\end{split}
\end{equation}
when
$$(\alpha \otimes 1 + 1 \otimes \tau(\alpha))\big( r_{\h_0,\h_0} +
r_{\l,\h_0}-r_{\l,\h_0}^{21} \big)=\frac{1}{2}(\alpha + \tau(\alpha))
\Omega_{\h}.$$
In particular, $\tilde{r}(\lambda)$ interpolates the constant r-matrix
obtained from the Belavin-Drinfeld triple $(\Gamma_{1}=\alpha,
\Gamma_{2}=\beta, \tau: \alpha \mapsto \beta)$ at $(\alpha,\lambda) \to
\infty$ and the r-matrix obtained from $(\Gamma_1=\beta, \Gamma_2=\alpha,
\tau: \beta \mapsto \alpha)$ at $(\alpha, \lambda) \to - \infty$.
\paragraph{Remark:} The generalization of this example to
$\g=\mathfrak{sl}_{2n+1}$ is the following. Fix a polarization and let
$\l=\C h_{\rho}$. Denote by $\Delta$ the root system and by $\Pi=(\alpha_1,
\ldots \alpha_{2n})$ the set of positive simple roots. Let $i: \alpha_{k}
\mapsto \alpha_{2n+1-k}$ be the involution of the Dynkin diagram. The
dynamical r-matrix obtained from the generalized Belavin-Drinfeld triple
$(\Gamma_{1}=\Gamma_{2}=\Pi,\tau=i)$ interpolates the constant r-matrices
obtained from the Belavin-Drinfeld triples $(\Gamma_1=(\alpha_1, \ldots
\alpha_n), \Gamma_2=(\alpha_{n+1}, \ldots \alpha_{2n}),\, \tau=i)$ and
$(\Gamma_1=(\alpha_{n+1}, \ldots \alpha_{2n}), \Gamma_2=(\alpha_{1}, \ldots
\alpha_{n}),\, \tau=i^{-1})$.
\subsection{Permutation dynamical r-matrices}
\paragraph{}Consider $\g=\mathfrak{sl}_{2n}$, and let $\Pi=(\alpha_{1},
\ldots \alpha_{2n-1})$ denote a system of simple roots. For any $\sigma \in
S_n$, we can construct a generalized Belavin-Drinfeld triple by setting
$\Gamma_{1}=\Gamma_{2}=(\alpha_1, \alpha_3, \ldots \alpha_{2n-1})$ and
$\tau: \alpha_{2k-1} \mapsto \alpha_{2 \sigma (k) -1}$.
\paragraph{Acknowledgements:}I heartily thank Pavel Etingof for his
encouragements, constant help and for his communicative enthousiasm for mathematics. I also
thank A. Varchenko and P. Etingof for sharing their work with me before
publication, and Hung Yean Loke and Vadik Vologodsky for interesting 
discussions.
\small{}
\end{document}